\documentclass[11pt]{article}

\usepackage{amsthm,amsmath,amssymb}
\usepackage[hypertex]{hyperref}
\usepackage[latin1]{inputenc}
\usepackage[english]{babel}
\usepackage[mathscr]{euscript}
\usepackage{mathscinet}
\usepackage[all]{xy}

\newtheorem{theo}{Theorem}[section]
\newtheorem{prop}[theo]{Proposition}
\newtheorem{lemma}[theo]{Lemma}
\newtheorem{cor}[theo]{Corollary}
\newtheorem{defi}[theo]{Definition}
\newtheorem{remark}[theo]{Remark}
\newtheorem{ex}[theo]{Example}
\newtheorem{conjet}[theo]{Conjecture}

\def\K{\mathbb{K}}
\def\u{\mathbf{u}}
\def\v{\mathbf{v}}
\def\X{\chi}
\def\l{\ell}
\def\bz{\mathbb{Z}}
\def\bn{\mathbb{N}}

\DeclareMathOperator{\ann}{Ann}
\DeclareMathOperator{\im}{im}

\title{\textbf{Graded Betti Numbers of the Logarithmic Derivation Module}}
\author{
{\Large M.~A.~Marco-Buzunáriz\footnote{Partially supported by ERC Starting Grant TGASS,
MTM2007-67908-C02-01 and ``E15 Grupo Consolidado Geometr\'ia''.}} \and
{\Large J.~Martín-Morales\footnote{Partially supported by MTM2007-67908-C02-01,
FQM-333 and ``E15 Grupo Consolidado Geometr\'ia''.}}
}

\date{\emph{
\begin{scriptsize}
\begin{tabular}{ccc}
ICMAT && Department of Mathematics-I.U.M.A.\\
 CSIC-Complutense-Aut\'onoma-Carlos III && University of Zaragoza (Spain)\\
mmarco@unizar.es && jorge@unizar.es
\end{tabular}
\end{scriptsize}}}

\begin{document}

\renewcommand{\labelenumi}{\bf \theenumi)}
\maketitle
\thispagestyle{empty}

\begin{abstract}
Let $Q\in \K[x_1,\ldots,x_n] = S$ be a homogeneous polynomial of degree~$d$. The freeness of the logarithmic derivation module, $D(Q)$, and of its natural generalizations, has been widely studied.
In the free case, $D(Q) \simeq \bigoplus_{i=1}^n S(-d_i)$ where the 
$d_i$'s are the exponents of the module; and as a direct consequence of the Saito-Ziegler criterion, the formula $d = \sum_i d_i$ holds. In this paper we give a generalization of this formula in the non-free case. Moreover, we show that an equivalent formula is also true in the quasi-homogeneous case, and show to what extent it can be generalized for arbitrary polynomials.

\end{abstract}

\section*{Introduction}

Let $\K$ be a field of characteristic zero. Let us denote by $S: = \K[x_1,\ldots,x_n]$
the ring of polynomials in $n$ variables with coefficients in $\K$. This ring is also an infinite dimensional vector space over $\K$. Let $Der_\K (S):=\{\theta :S\rightarrow S\mid \theta(ab)=a\theta(b)+b\theta(a)\}$ be the module of the $\K$-linear derivations of $S$. It can be shown that $Der_\K(S)$ is a free $S$-module with basis $\{\partial_1,\ldots,
\partial_n\}$, where $\partial_i (g):= \frac{\partial g}{\partial x_i}$. In other words,
every derivation in $S$ can be written in a unique way in the form
$\delta = \sum_{i=1}^n a_i \partial_i$ with $a_i\in S$.

Given a polynomial $f\in S$, a derivation $\delta\in Der_\K(S)$ is said to be a {\em logarithmic derivation} with respect to $f$ if $\delta(f) = h f$ for some $h\in S$. The set of logarithmic derivations with respect to $f$ is denoted by
$Der(-log\ f)$. It is easy to check that $Der(-log\ f)$ is a submodule of $Der_\K(S)$. If $f = gh$, then one has that
$$
  Der(-log\ f) = Der(-log\ g) \cap Der(-log\ h),
$$
even if $g$ and $h$ do have common factors.

The previous definition does not distinguish between a polynomial and its reduced. For this reason, it is convenient to define a slightly different module as follows: let $f=f_1^{e_1} f_2^{e_2}\cdots f_r^{e_r}$ with $f_i$ irreducible for each $i$. Consider 
$$
  D(f) := \bigcap_{i=1}^r D(f_i;\, e_i),
$$
where $D(g;\, k) = \{ \delta\in Der_\K(S) \mid \exists h\in S,\ \delta(g) = h g^k \}$. This definition coincides with  the logarithmic derivation module when $f$ is reduced, since $D(g,1): = Der(-log\ g)$. We will call this module the \emph{(generalized) logarithmic derivation module} of $f$.

It is easy to prove that, when $D(f)$ is free, its rank is $n$. The so-called \emph{Saito criterion} characterizes the bases of the module $Der(-log f)$ when free.

\begin{theo}[Saito criterion {\cite[Th.~4.19]{Saito80}}]
Let $f\in S$ be an arbitrary polynomial, and let $\delta_1,\ldots,\delta_n$ be
derivations in $D(f)$. The following statements are equivalent:
\begin{enumerate}
  \item $\{\delta_1,\ldots,\delta_n\}$ is a basis of $D(f)$.
  \item $\det [\delta_1|\cdots|\delta_n] = c\cdot f$, for some $c\in\K^{*}$.
\end{enumerate}
\end{theo}

Ziegler studied the module $D(Q)$ for arrangements and gave a generalization
of the previous result in \cite[p.~351]{Ziegler89}.

The case when the module is free has been widely studied; particularly in the setting of arrangements (see \cite{OT92},\cite{Saito80},\cite{Schenck09} among others) and in relation with the logarithmic comparison theorem  (a survey can be found in \cite{Torrelli07}).

When the polynomial $Q$ is homogeneous of degree $d$ the Saito-Ziegler criterion allows us to decide algorithmically whether
$D(Q)$ is free as $S$-module. In that case $D(Q) \simeq \bigoplus_{i=1}^n S(-d_i)$
where the numbers $d_i$ are called the {\em exponents} of the module
and are precisely the degrees of a homogeneous basis. As a consequence, we have that $d = \sum_{i=1}^n d_i$.

In this paper we will study the relationship between $d$ and the graded betti numbers when the module is not necessarily free,
giving a general formula that particularizes to the previous one in the free case.
More precisely, the following result will be proved (see notation in Section~\ref{homogDQ}).

\begin{theo}\label{mainTheo}
Let $Q\in S^\u$ be a $\u$-homogeneous polynomial with $\u\in\mathbb{Z}^n_{> 0}$, and $\v\in\mathbb{Z}^n$
a vector such that $\u+\v = (k,\ldots,k) =: k\cdot{\bf 1}$. Then $D^{(\u,\v)}(Q)$
is a graded $S^\u$-submodule of $D^{(\u,\v)} = Der_\K^{(\u,\v)}(S^\u)$. Moreover, if
$$
  0\longleftarrow D^{(\u,\v)} (Q)\longleftarrow \bigoplus_{i=1}^{r_0} S^\u(-d_i^0) \longleftarrow
  \cdots \longleftarrow \bigoplus_{i=1}^{r_\l} S^\u(-d_i^\l) \longleftarrow 0
$$
is a $(\u,\v)$-homogeneous free resolution of $D^{(\u,\v)} (Q)$, then
$\deg^\u(Q)+|\v|$ is the alternating sum of the exponents that appear in the resolution.
That is, $\deg^\u(Q)+|\v| = \sum_{p=0}^\l (-1)^{p} \sum_{i=1}^{r_p} d_i^p$. 
\end{theo}
In particular, taking a minimal free resolution, one has that
$$
  \deg^\u(Q) + |\v| = \sum_{p=0}^\l (-1)^p \sum_{j\in\mathbb{Z}} (j+p)\, b_{j,p},
$$
where $b_{j,p} = b_{j,p}(D^{(\u,\v)}(Q))$ are the betti numbers of $D^{(\u,\v)} (Q)$
as $(\u,\v)$-graded module.

The key step in the proof of this theorem is the interpretation of the exponents that appear in a
homogeneous free resolution of a graded module in terms of its Hilbert-Poincar\'e series.

The paper is organized as follows. In Section \ref{minRes}, we recall some notions of graded modules,
homogeneous free resolutions and graded betti numbers. Section \ref{HPseries} contains the definition of the Hilbert-Poincar\'e series and some examples that ilustrate how to compute it for the modules we are studying. In section \ref{homogDQ} we stablish conditions for the module $D(Q)$ to be quasi-homogeneous (and equivalenty, for the betti numbers to be well defined). The relationship between the Hilbert-Poincar\'e series and the betti numbers is discused in Section \ref{interpreHPS}. This relationship is used in Section \ref{computeXDQ} to prove the main result. Finally, in Section  \ref{generalCase} we show to what extent this result can be generalized for the non-graded case.

\section{Minimal Resolutions and Betti Numbers}\label{minRes}%

The following two sections contain the definitions of some classical objects that will be used later. Details can be
found, for instance, in \cite{Eisenbud95} or \cite{GP02}.

\begin{defi}
Let $\K$ be a field, and $R$ a $\K$-algebra. A \textbf{graded structure} on $R$ is a decomposition as $\K$-vector space $R=\bigoplus_{n\in\bn}R_n$
such that $R_i\cdot R_j\subseteq R_{i+j}$. An algebra is said to be \textbf{graded} when it admits a graded
structure. In that case, the subspaces $R_i$ are called the \textbf{homogeneous parts} of $R$.
\end{defi}

In particular, the following two cases of graded algebras will be used.

\begin{ex}\label{ejemHom}
The ring of polynomials $S=\K[x_1,\ldots,x_n]$ is graded, being its $i$-th homogeneous part precisely the space of the homogeneous polynomials of degree $i$.
\end{ex}

\begin{ex}\label{ejmmodqh}
 The same ring $S=\K[x_1,\ldots,x_n]$ admits different graded structures by choosing a weight vector $\u=(u_1,\ldots,u_n)\in\bz_{>0}^n$, and choosing as $S_i$ the vector space spanned by the monomials $x_1^{a_1}\cdots x_n^{a_n}$ such that $a_1u_1+\cdots+a_nu_n=i$. This structure reflects the idea that the variable $x_i$ has degree $u_i$. In this case, the polynomials that lie in $S_i$ are called \textbf{quasi-homogeneous} with respect to $\u$. We will use the notation $S^\u$ to emphasize that we are using this graded structure if necesary.
\end{ex}

Let us recall the definition of graded module.
\begin{defi}
 Let $R=\bigoplus_{n\in\bn}R_n$ be a graded $\K$-algebra, and let $M$ be a $R$-module. A \textbf{graded structure} on $M$ is a decomposition as $\K$-vector space $
M=\bigoplus_{i\in\bz}M_i$ such that $R_i\cdot M_j\subseteq M_{i+j}$. As above, we call a module with a graded structure a \textbf{graded module}, and the $M_i$ will be called the \textbf{homogeneous parts} of $M$. 
\end{defi}

A submodule $N$ of a graded module $M$ inherits the graded structure if the homogeneous parts of each element of $N$ are also in $N$. In that case, the quotient $M/N$ also inherits the graded structure.

\begin{defi}
Given two graded modules, $M$ and $N$, a morphism $\varphi:M\rightarrow N$ is said to be \textbf{homogeneous of degree $d$} if $\varphi(M_i)\subseteq N_{i+d}$.
\end{defi}

\begin{ex}
 Given a graded module $M$ and an integer number $d$, one can define a new graded structure (denoted by $M(d)$) just by shifting the degrees of $M$ by $d$, that is $M(d)_i:=M_{i+d}$. In particular, if we consider $S$ as module over itself, the identity map $S\rightarrow S(-d)$ is a homogeneous morphism of degree $d$.
\end{ex}

\begin{defi}
Let $R$ be a graded $\K$-algebra and let $M$ be a finitely generated graded $R$-module. A \textbf{homogeneous free resolution} of $M$
is a resolution
$$
  0 \longleftarrow M \longleftarrow F_0 \stackrel{\varphi_1}{\longleftarrow} F_1 \longleftarrow \cdots
  \longleftarrow F_{k-1} \stackrel{\varphi_k}{\longleftarrow} F_{k} \longleftarrow \cdots
$$
such that each $F_k$ is a finitely generated free $R$-module,
$$
  F_k = \bigoplus_{j\in\mathbb{Z}} R(-j-k)^{b_{j,k}},
$$
and each $\varphi_k$ is a homogeneous morphism of degree $0$.

Such a resolution is said to be \textbf{minimal} if $\varphi_k(F_k) \subset
\mathfrak{m} F_{k-1}$, being $\mathfrak{m}$ the ideal of $R$ generated by homogeneous parts of positive degree. 
In this situation, the numbers 
$b_{j,k}(M)$ are called the \textbf{graded betti numbers} of $M$
and $b_k(M) = \sum_j b_{j,k}(M) = rank(F_k)$ is called \textbf{$k$-th
betti number} of $M$.
\end{defi}

A minimal resolution always exists,
and the graded betti numbers are independent of the chosen minimal homogeneous free resolution.

\section{Hilbert-Poincaré Series}\label{HPseries}%

\begin{defi}
 Let $R$ be a graded $\K$-algebra and let $M$ be a graded $R$-module. The formal series in $t$
$$\sum_{i\in\bz}(dim_{\K}M_i)\cdot t^i
$$
is called the \textbf{Hilbert-Poincaré series} of $M$, and will be denoted by $HP_M(t)$.
\end{defi}

\begin{ex}
 If $R=\K[x]$ with the natural grading and we consider $R$ as a module over itself, we have that $R_i=\K\cdot x^i$. In this case, $$HP_{\K[x]}(t)=1+t+t^2+\cdots=\frac{1}{1-t}.$$ Analogously, it can be proved that $$HP_{\K[x_1,\ldots,x_n]}(t)=\frac{1}{(1-t)^n}.$$
\end{ex}

The previous formula can be generalized to the quasi-homogeneous case as follows.

\begin{ex}\label{ejmHPqh}
If $S^\u=\K[x_1,\ldots,x_n]$ is endowed with the graded structure of Example \ref{ejmmodqh} (with a weight vector $\u=(u_1,\ldots,u_n)\in\mathbb{Z}^n_{>0}$), it is readily that $$HP_{S^\u}(t)=\frac{1}{\prod_{i=1}^n(1-t^{u_i})}.$$
\end{ex}

Note that $t^{-d} HP_M(t)=HP_{M(d)}(t)$. This fact, together with the well-known fact that the Hilbert-Poincaré series is an additive functor, allows us to compute it for every module by looking at a graded free resolution.

The Hilbert-Poincaré Series is related to the exponents as follows.
\begin{remark}\label{ejemHPlibre}
Let $M = \bigoplus_{i=1}^r S^\u(-d_i)$ be a free graded module.
 Then $$HP_M (t) = \frac{\sum_{i=1}^r t^{d_i}}{\prod_i (1-t^{u_i})}.$$
In particular $\big[\prod_i(1-t^{u_i}) HP_M (t)\big]'(1) = \sum_{i=1}^r d_i$.
\end{remark}


The following result is due to Hilbert and can be found in the standard literature (cf. \cite{Eisenbud95}, \cite{GP02}).
\begin{lemma}
 Let $S^\u$ be as in Example~\ref{ejmHPqh}, and $I$ a $\u$-homogeneous ideal of $S^\u$. The order of the pole of $HP_{S^\u/I}(t)$ at $t=1$ coincides with the Krull dimension of the ring $S^\u/I$.
\end{lemma}

\section{Homogeneity of $D(Q)$}\label{homogDQ}

Throughout this section $Q\in S$ will be a quasi-homogeneous polynomial with respect to a vector $\u\in\mathbb{Z}^n_{>0}$ as in Example~\ref{ejmmodqh}.
We can give a graded structure on $Der_\K(S)$ in such a way that $D(Q)$ is a graded $S^\u$-submodule as follows. The degree of a derivation is given by $\deg(a_i\partial_i) =
\deg^\u(a_i)+v_i$ where $a_i\in S^\u$. That is, the partial derivatives $\partial_i$ have degree $v_i$, and hence $Der_\K(S)\simeq \bigoplus_{i=1}^n S^\u(-v_i)$ as a graded $S^u$-module. We will denote this graded structure by $Der_\K^{(\u,\v)}(S^\u)$, or simply $D^{(\u,\v)}$, where $\v=(v_1,\ldots,v_n)$. If the graded structure is not relevant, the vectors $\u,\v$ will be omited.

In what follows we will see for which values of $\u$ and $\v$ the module $D(Q)$ turns out to be a graded $S^\u$-submodule of $D^{(\u,\v)}$.

\begin{ex}
Let $Q=x^2+y^2\in\K[x,y] = S$. It is $\u$-homogeneous with $\u=(u,u)$. Consider the Euler derivation $\mathcal{E} = x\partial_x + y\partial_y$. Clearly $\mathcal{E}\in D(Q)$, so if $D(Q)$ is a graded submodule of $D^{(\u,\v)}$, the homogeneous parts of $\mathcal{E}$ must also be in $D(Q)$.
Since $x\partial_x,\, y\partial_y\notin D(Q)$, the only way $D(Q)$ can be graded is if
the Euler derivation is homogeneous, and hence the vector $\v$
must have the form $\v=(v,v)$.
\end{ex}


\begin{ex}
Let $Q=x^{e_1} y^{e_2}\in\K[x,y] = S^\u$ for some $\u\in\bz_{>0}^2$ and consider the derivations
$\delta_1 = x^{e_1}\partial_x$,
$\delta_2 = y^{e_2}\partial_y$. Clearly
$\delta_1\,$ and $\,\delta_2$ are derivations that lie in $D(Q)$.
By Saito's criterion, we can see that, in fact, $\{\delta_1,\delta_2\}$
form a basis of $D(Q)$. Now for every $\v=(v_1,v_2)\in\mathbb{Z}^2$
the module $D^{(\u,\v)} (Q)$ is a graded submodule of $Der_\K^{(\u,\v)}(S^\u)$.
\end{ex}

\begin{lemma}\label{lemmaQH1}
For all $f\in S$, one has $\ann_S(D/D(f)) = S\langle f \rangle$.
\end{lemma}

\begin{proof}
One of the inclusions is obvious since $f\cdot D \subseteq D(f)$.
Let us consider $f=f_1^{e_1}\cdots f_r^{e_r}$ the
decomposition of $f$ into irreducible factors and
let $g\in \ann_S (D/D(f))$. Then $g\cdot D\subseteq D(f)$
and thus $f_i^{e_i}$ divides $g\delta(f_i)$ for all $\delta\in D$ and
for all $i$. To finish the proof it is enough to show that there exists
$\delta\in D$ such that $\delta(f_i)$ and $f_i^{e_i}$ do not have
common components for all $i$.

If $f$ is $x_1$-general then so is $f_i$ for all $i$ and $\delta
= \frac{\partial}{\partial x_1} = \partial_{x_1}$ can be chosen.
Otherwise there exists $\varphi: S \to S$
a change of coordinate such that $\varphi(f) = g$ where $g$ is now $x_1$-general.
Let us denote $g_i = \varphi(f_i)$ and $\psi = \varphi^{-1}$. Thus the polynomials
$\frac{\partial g_i}{\partial x_1}$ and ${g}_i^{e_i}$ do not have common factors and
therefore $\psi(\frac{\partial g_i}{\partial x_1})$ and $\psi({g}_i^{e_i})= f_i^{e_i}$
do not either. Note that
$$
  \psi\left(\frac{\partial g_i}{\partial x_1}\right) =
  \sum_{k=1}^n \psi\!\left(\frac{\partial \varphi(x_k)}{\partial x_1}\right)\cdot
  \frac{\partial f_i}{\partial x_k}.
$$
Hence $\delta = \sum_{k} \psi\Big(\frac{\partial \varphi(x_k)}{\partial x_1}\Big) \partial_{x_k}$
verifies the aforementioned condition and the claim follows.
\end{proof}

Recall that the dimension of a module is related to the dimension of its annihilator.
Specifically, the above result shows that the logarithmic derivation module $D(Q)$ always has
dimension $n-1$.

\begin{lemma}\label{lemmaQH2}
Let $Q\in (S^\u)_d$ be a homogeneous polynomial of degree $d$ and $\delta\in (D^{(\u,\v)})_j$
a homogeneous derivation of degree $j$ with $\u+\v = k\cdot {\bf 1}$. Then
$\delta(Q) \in (S^\u)_{d+j-k}$, i.e. it is homogeneous of degree $d+j-k$.
\end{lemma}

\begin{proof}
Assume $Q = x_1^{\alpha_1}\cdots x_n^{\alpha_n} = x^\alpha$ is a monomial
with $\alpha\cdot \u :=\sum_i \alpha_i u_i= d$. The general case follows obviously from this one.
The derivation $\delta$ can be written in the form $\delta = \sum_i h_i \partial_{x_i}$
where $h_i\in (S^\u)_{j-v_i}$ has degree $j-v_i$. The degree of $h_i x^{\alpha-e_i}$
in $S^\u$ is $j-v_i + (\alpha-e_i)\cdot \u = j-v_i + d - u_i = d+j-k$.
Hence $\delta(Q) = \sum_i \alpha_i h_i x^{\alpha-e_i} \in (S^\u)_{d+j-k}$.
\end{proof}

\begin{prop}\label{equivcuasihom}
\label{prophomog}
Let $Q\in S$ and $\u\in\mathbb{Z}^n$.
The following conditions are equivalent:
\begin{enumerate}
\item The polynomial $Q$ is quasi-homogeneous with respect to the weight vector~$\u$.
\item The module $D(Q)$ is a graded $S^\u$-submodule of $D^{(\u,\v)}$ for
all $\v\in\mathbb{Z}^n$ with $\u+\v = k\cdot {\bf 1}$.
That is $D(Q)$ is homogeneous with respect to the weight vector $(\u,\v)$.
\end{enumerate}
\end{prop}

\begin{proof}
Assume $Q$ is $\u$-homogeneous of degree $d$ and consider $\delta \in D(Q) \subseteq D^{(\u,\v)}$.
Let $\delta = \sum_j \delta_j$ be the decomposition of $\delta$ into homogeneous
parts, that is $\delta_j \in (D^{(\u,\v)})_j$. Let us consider $Q=Q_1^{e_1}\cdots Q_r^{e_r}$
the decomposition of $Q$ into irreducible factors and denote by $d_i$ the degree
of $Q_i$. Then for each $i$ there exists $h_i$
with $\delta(Q_i) = h_i Q_i^{e_i}$. Every $h_i$ can be written as $h_i = \sum_j h_{ij}$
where $h_{ij}$ has degree $j$. Taking the homogeneous parts in the expression
$$
  \sum_j h_{ij} Q_i^{e_i} = h_i Q_i^{e_i} = \delta(Q_i) = \sum_j \delta_j (Q_i)\ \in\ S^\u,
$$
we have that $\delta_j(Q_i) = h_{i,j-d_i (e_i-1)-k}\cdot Q_i^{e_i}$.
In particular $\delta_j\in D(Q)$ for all $j$ and hence $D(Q)$ is $(\u,\v)$-homogeneous.

Assume now that the condition (2) is satisfied.
The $\u$-homogeneity of $Q$ can be easily deduced from Lemma \ref{lemmaQH2},
since the annihilator of a graded module is always a homogeneous ideal with
the inherited graduation.
\end{proof}

When these conditions hold, we will use the notation $D^{(\u,\v)}(Q)$ to emphasize the graded structure.

\begin{remark}Note that $D(Q)$ can be computed using Gr\"obner bases in the ring
of polynomials $S$. If $Q$ is quasi-homogeneous, all operations required in the computation of $D(Q)$ preserve the quasi-homogeneity.
This actually provides another proof for the first implication of Proposition~\ref{equivcuasihom}.
\end{remark}

\section{Interpretation in terms of Hilbert-Poincaré series}\label{interpreHPS}

Given a graded submodule 
$M\subseteq D^{(u,v)}=Der_\K^{(\u,\v)} (S^\u)$, consider $HP_M (t)$ the Hilbert-Poincaré series
of $M$. Now define $\X(M):= \big[\prod(1-t^{u_i}) HP_M (t)\big]'(1)\in\mathbb{Z}$.
The following properties of $\X$ are a direct consequence of the properties verified by 
the Hilbert-Poincaré series.

\begin{lemma}
\label{propchi}
The invariant $\X$ has the following properties:
\begin{itemize}
  \item $\X$ is additive,
  \item $\X(S^{\u}(-d)) = d$,
  \item $\X(D^{(\u,\v)})=|\v| := v_1 + \cdots + v_n$.
\end{itemize}
\end{lemma}

From now on, $\u=(u_1,\ldots,u_n)\in \mathbb{Z}^n_{>0}$ and $\v = (v_1,\ldots,v_n)\in\mathbb{Z}^n$ will be vectors such that
$D^{(\u,\v)} (Q)$ is a graded $S^\u$-submodule of $D^{(\u,\v)}$.
The purpose of this section is to reduce Theorem~\ref{mainTheo}
to the computation of $\X(D^{(\u,\v)} (Q))$. The following two results achieve this.

\begin{prop}\label{simplificacion}
Let $Q\in S^\u$ be a quasi-homogeneous polynomial and consider $D^{(\u,\v)} (Q)\subseteq D^{(\u,\v)}$
its module of logarithmic derivations. Suppose that
\begin{equation}\label{res1}
  0\longleftarrow D^{(\u,\v)} (Q)\longleftarrow \bigoplus_{i=1}^{r_0} S(-d_i^0) \longleftarrow
  \cdots \longleftarrow \bigoplus_{i=1}^{r_\l} S(-d_i^\l) \longleftarrow 0
\end{equation}
is a free homogeneous resolution of $D^{(\u,\v)} (Q)$. Then
$\X(D^{(\u,\v)} (Q))$ is the alternating sum of the exponents that appear in (\ref{res1}).
That is, $\X(D^{(\u,\v)} (Q)) = \sum_{p=1}^\l (-1)^{p} \sum_{i=1}^{r_p} d_i^p$.
\end{prop}

\begin{proof}
Take $F_{p} = \bigoplus_{i=1}^{r_p} S(-d_i^{p})$ for $0\leq p\leq \l$.
Since $\X$ is additive, it is readily seen that
$\X (D(Q)) = \sum_{p=0}^l (-1)^p\, \X (F_p)$.
On the other hand, $$\X(F_p) = \X\big(\bigoplus_{i=1}^{r_p} S(-d_i^p)\big)
= \sum_{i=1}^{r_p} \X\big(S(-d_i^p)\big) = \sum_{i=1}^{r_p} d_i^p$$ and
the claim follows.
\end{proof}

\begin{cor}\label{cor1}
Let $b_{j,p} = b_{j,p}(D(Q))$ be the graded betti numbers of $D^{(\u,\v)}(Q)$ with respect to the grading given
by $(\u,\v)$. Then
$$
  \X(D^{(\u,\v)} (Q)) = \sum_{p=0}^\l (-1)^p \sum_{j\in\mathbb{Z}} (j+p)\, b_{j,p},
$$
\end{cor}

\begin{proof}
If the resolution of Proposition~\ref{simplificacion} is minimal, then, by
definition of betti numbers,
$F_p = \bigoplus_{i=1}^{r_p} S(-d_i^p) =
\bigoplus_{j\in\mathbb{Z}} S(-j-p)^{b_{j,p}}$.
Applying $\X$ to the previous equality one obtains that $\sum_{i} d_i^p =
\sum_{j} (j+p)\, b_{j,p}$.
\end{proof}

In the next section we will see that $\X(D^{(\u,\v)} (Q)) = \deg(Q) + |\v|$.

\section{Computation of the invariant $\X(D^{(\u,\v)}(Q))$}\label{computeXDQ}

The following lemma is key for the computation of $\X(D^{(\u,\v)} (Q))$.

\begin{lemma}\label{clave}
Let $I\subseteq S^\u$ be a homogeneous ideal. If $I$ contains two polynomials
with no common factors, then there exists
 $H(t)\in \mathbb{Z}[t]$ such that $\prod_i(1-t^{u_i}) HP_{S^\u/I}(t) = (1-t)^2 H(t)$.
In particular $\X((S^\u/I)(-d))=0$.
\end{lemma}

\begin{proof}
The Hilbert-Poincaré series of $S^\u/I$ can always be written as
 $HP_{S^\u/I}(t) = (1-t)^{(n-s)}\cdot G(t)/\prod(1-t^{u_i})$ where
$G(1)\neq 0$ and $s=\dim(S^\u/I)$.
We want to prove that $n-s\geq 2$. Since $\sqrt{I}$ contains two polynomials
without common components $n-2\geq \dim(S^\u/\sqrt{I}) = \dim(S^\u/I) = s$.

For the last statement, note that $$\prod_i (1-t^{u_i}) HP_{(S^\u/I)(-d)}(t) = t^d (1-t)^2 H(t),$$
and $[t^d (1-t)^2 H(t)]'(1)=0$.
\end{proof}

%
%

\begin{lemma}\label{construccionMi}
Let $M$ be a graded submodule of $D^{(\u,\v)}$ and consider
$M_i = M + S^\u\langle \partial_1,\ldots,\partial_i\rangle$
for $i=0,\ldots,n$. ($M_0=M,\ M_n = D^{(\u,\v)}$). Then
$$
  \X(D^{(\u,\v)}/M) = \sum_{i=1}^n \X (M_i/M_{i-1}).
$$
\end{lemma}

\begin{proof}
Recall that, given three modules, $N\subset M \subset L$
we can define the short exact sequence $\, 0\to M/N \to L/N \to L/M \to 0\,$
given by inclusion and projection respectively.
Since $M_{i-1} \subset M_i \subset D^{(\u,\v)}$ we get that
$$
  0 \longrightarrow M_i/M_{i-1} \longrightarrow D^{(\u,\v)}/M_{i-1} \longrightarrow
  D^{(\u,\v)}/M_i \longrightarrow 0,\quad i=1,\ldots,n.
$$

Moreover, all the morphisms are homogeneous of degree $0$.
Now we can use these exact sequences and the additivity of $\X$
to obtain the required expression (recall that $M_0 = M$, $M_n = D^{(\u,\v)}$).
\begin{equation*}
\begin{split}
  \X(D^{(\u,\v)}/M) & = \X(M_1/M_0) + \X(D^{(\u,\v)}/M_1) =\\
  & = \X(M_1/M_0) + \X(M_2/M_1) + \X(D^{(\u,\v)}/M_2) =\\
  & \hspace{0.175cm} \vdots\\
  & = \X(M_1/M_0) + \cdots + \X(M_n/M_{n-1}) + \X(D^{(\u,\v)}/M_{n}) =\\
  & = \sum_{i=1}^n \X(M_i/M_{i-1})
\end{split}
\end{equation*}
\end{proof}

\begin{prop}\label{noCompComunes}
Let $Q_1$ and $Q_2\in S$ be two quasi-homogeneous polynomials. If $Q_1^{red}$ and
$Q_2^{red}$ have no common factors, then $\X(D^{(\u,\v)}(Q_1)+D^{(\u,\v)}(Q_2))=|\v|$.
\end{prop}

\begin{proof}
Let $M = D^{(\u,\v)}(Q_1) + D^{(\u,\v)}(Q_2)$; consider the modules
$M_i = M + S^\u\langle \partial_1,\ldots,\partial_i\rangle$, $i=0,\ldots,n$.
We will prove that $\X(M_i/M_{i-1}) = 0$, $i=1,\ldots,n$. Using Lemma~\ref{construccionMi} 
and the additivity of $\X$ we get that
$\X(M) = \X(D^{(\u,\v)}) = \sum \X(S^\u(-v_i)) = |\v|$.

Now observe that
$$
  \frac{M_i}{M_{i-1}} = \frac{M_{i-1} + S^\u\langle \partial_i \rangle}{M_{i-1}}
  \cong \frac{S^\u\langle \partial_i \rangle}{S^\u\langle \partial_i\rangle \cap M_{i-1}}
  \cong \frac{S^\u}{I_i}(-v_i),\quad i=1,\ldots,n.
$$
The first isomorphism comes from the (Second) Isomorphism Theorem, and the second one is induced by the identification
$S^\u\langle \partial_i \rangle\cong S^\u(-v_i)$. Particularly $I_i$ is the image of
$S^\u\langle \partial_i \rangle\cap M_{i-1}$ by the morphism
$S^\u\langle \partial_i \rangle \to S^\u(-v_i)$ given by $h\partial_x \mapsto h$
for $h\in S^\u$. These two morphisms are homogeneous of degree 0 because
$\partial_i$ has degree $v_i$.

The derivation $Q_1\partial_i$ is clearly in $D^{(\u,\v)}(Q_1)\subseteq M \subseteq M_{i-1}$.
Therefore $Q_1\partial_i\in S^\u\langle \partial_i\rangle \cap M_{i-1}$ and hence $Q_1\in I_i$.
For the same reason $Q_2$ is also in $I_i$. We have just proved that
$I_i$ is a homogeneous ideal that contains two polynomials without common factors. By Lemma~\ref{clave} 
and the previous isomorphisms
$\X(M_i/M_{i-1}) = \X((S^\u/I_i)(-v_i)) = 0$.
\end{proof}

\begin{theo}\label{XDQhom}
$\X(D^{(\u,\v)}(Q)) = \deg^\u(Q)+|\v|$.
\end{theo}

\begin{proof}
If $Q$ is a constant polynomial, then $D^{(\u,\v)}(Q)=D^{(\u,\v)}$ and the result is a direct consequence of Lemma~\ref{propchi}.

Now assume that $Q$ has a positive degree, and let $Q=Q_1^{e_1}\cdots Q_r^{e_r}$
its irreducible factor decomposition. We shall complete the proof by induction on $r$.

\vspace{15pt}

\framebox{$r=1$}\,. Let $M=D^{(\u,\v)}(Q)=D^{(\u,\v)}(Q_1^{e_1})$ and consider the modules
$M_i = M + S^\u\langle \partial_1,\ldots,\partial_i \rangle$,
$i=0,\ldots,n$. We will prove that $\X(M_i/M_{i-1})=0$, for $i=2,\ldots,n$
and that $\X(M_1/M) = -d$. By Lemma \ref{construccionMi}, this would prove that  $\X(D^{(\u,\v)}/M) = -d$ and, by Lemma~\ref{propchi}, $\X(M)=d+|\v|$.

Without loss of generality, we may assume that $Q_1$ is $x_1$-general of order $deg(Q_1)$.

Take $i\in\{2,\ldots,n\}$. As in the proof of Proposition
\ref{noCompComunes} we have the following graded isomorphisms:
$$
  \frac{M_i}{M_{i-1}} = \frac{M_{i-1} + S^\u\langle \partial_i \rangle}{M_{i-1}}
  \cong \frac{S^\u\langle \partial_i \rangle}{S^\u\langle \partial_i\rangle \cap M_{i-1}}
  \cong \frac{S^\u}{I_i}(-v_i),\quad i=1,\ldots,n,
$$
where $I_i$ is the image of $S^\u\langle \partial_i\rangle \cap M_{i-1}$
by the graded morphism $S^\u\langle \partial_i \rangle \to S^\u(-v_i)$ given by
$h\partial_i \to h$. It is clear that $Q_1^{e_1}\partial_i\in D^{(\u,\v)}(Q_1^{e_1}) = M
\subseteq M_{i-1}$. This implies $Q_1^{e_1}\partial_i\in
S^\u\langle \partial_i \rangle \cap M_{i-1}$ and hence $Q_1^{e_1}\in I_i$.
On the other hand, the derivation $\frac{\partial Q_1}{\partial x_i}\partial_1 -
\frac{\partial Q_1}{\partial x_1}\partial_i \in M$. Since $\partial_1\in M_{i-1}$
one obtains that  $\frac{\partial Q_1}{\partial x_1}\partial_i\in M_{i-1}$
and thus $\frac{\partial Q_1}{\partial x_1} \in I_i$.
(Note that this last statement is not necessarily true for $i=1$, since in general
$\partial_1\notin M=M_0$). From the assumption that $Q_1$ is $x_1$-general
of order $deg(Q_1)$, one obtains that $Q_1^{e_1}$ and $\frac{\partial Q_1}{\partial x_1}$
have no common factors. By Lemma~\ref{clave} it follows that $\X(M_i/M_{i-1})) = \X((S^\u/I_i)(-v_i)) = 0$.

Let us now study the case  $i=1$, that is, we are interested in the module $M_1/M$. Analogously
 $M_1/M \cong S^\u\langle \partial_1 \rangle / S^\u\langle \partial_1 \rangle
\cap M$. Let us see that $S^\u\langle \partial_1 \rangle \cap M = S^\u\langle Q\partial_1\rangle$.
It is clear that $Q\partial_1\in M = D^{(\u,\v)}(Q)$. Let us suppose that $\delta = h\partial_1\in M$
with $h\in S^\u$. Then there exists $h_1\in S^\u$ such that $h\frac{\partial Q_1}{\partial x_1} =
h\partial_1(Q_1) = h_1 Q_1^{e_1}$. Since $Q_1^{e_1}$ and $\frac{\partial Q_1}{\partial x_1}$
have no common factors, it follows that
$Q_1^{e_1}$ divides $h$. So $\X(M_1/M) =
\X(S^\u\langle \partial_1 \rangle/S^\u\langle \partial_1 \rangle \cap M) =
\X(S^\u\langle \partial_1 \rangle)- \X(S^\u\langle \partial_1 \rangle \cap M) =\X(S^\u\langle \partial_1 \rangle)-\X(S^\u\langle Q\partial_1\rangle)=
v_1 - (v_1+d) = -d$.

\vspace{15pt}

\framebox{$<r \to r$}\,. Suppose $Q=Q_1^{e_1} Q_2$ where $Q_1$ y $Q_2$ have no common factors, and
by induction hypothesis, suppose the result is true for $Q_2$. 
Clearly $D^{(\u,\v)}(Q) = D^{(\u,\v)}(Q_1^{e_1})
\cap D^{(\u,\v)}(Q_2)$. The additivity of $\X$ implies
$$
  \X(D^{(\u,\v)}(Q))  = \X(D^{(\u,\v)}(Q_1^{e_1})) + \X(D^{(\u,\v)}(Q_2)) - \X(D^{(\u,\v)}(Q_1^{e_1}) + D^{(\u,\v)}(Q_2)).
$$
The first term of the sum is $\deg(Q_1^{e_1})+|\v|$ for the case $r=1$, the second one is $\deg(Q_2)+|\v|$ by induction hypothesis and the last one is $|\v|$ because of Lemma~\ref{noCompComunes}.
\end{proof}

\begin{remark}
One always has that
$\frac{\partial Q}{\partial x_i} \partial_{x_1} - \frac{\partial Q}{\partial x_1} \partial_{x_i} \in D(Q)$
and hence $\frac{\partial Q}{\partial x_1},Q\in I_i$, $i=2,\ldots,n$. Therefore when $Q$ is reduced,
since $D(Q)$ and $Der(-log\ Q)$ coincide, no induction is needed.
\end{remark}

We have proved the main result of this paper. Now assume that $Q$ is a $\u$-homogeneous polynomial and $D^{(\u,\v)}(Q)$ is a graded module. From Proposition~\ref{prophomog}, $D^{(\u,\v+{\bf{1}})}(Q)$ is also graded. Applying Theorem~\ref{XDQhom} to these two modules, we obtain:
$$
\X(D^{(\u,\v+{\bf{1}})}(Q))-\X(D^{(\u,\v)}(Q))=|\v+{\bf{1}}|-|\v|=n.
$$

Let us compute $\X(D^{(\u,\v)}(Q))$ using a free resolution $0\leftarrow D^{(\u,\v)}(Q) \leftarrow F_\bullet \leftarrow 0$ . On the other hand, the resolution $0\leftarrow D^{(\u,\v+{\bf 1})}(Q) \leftarrow F_\bullet(-1) \leftarrow 0$ can be used to compute $\X(D^{(\u,\v+{\bf{1}})}(Q))$. The difference between both resolutions is that the $S^\u(-i)$ of the first one appear as $S^u(-i-1)$ in the second one. So, when computing the alternating sum we obtain that $\sum_{p\geq 0} (-1)^p rk(F_p)=n$. That is, the alternating sum of the ranks of a free homogeneous resolution must be equal to the number of variables.
\begin{remark}
The previous statement can also be proved as follows:

Since $K(S)$ the fraction field of $S$ is a flat $S$-module,
 the rank of a finitely generated $S$-module $M$ defined as $rk(M) = K(S)\otimes_S M$ is an additive function.
If $M=D^{(\u,\v)}(Q)$ and we take a free homogeneous resolution
$$0\leftarrow M \leftarrow F_\bullet \leftarrow 0$$ then
$rk(M) = \sum_{p\geq 0} (-1)^p rk(F_p)$. Using the inclusion
$QD^{(\u,\v)} \subset D^{(\u,\v)}(Q) \subset D^{(\u,\v)}$ and the flatness of the $S$-module $K(S)$
 one can verify that $$n = rk(D^{(\u,\v)}(Q)) = \sum_{p\geq 0} (-1)^p rk(F_p).$$
\end{remark}

Note that when the weight vector $\u$ contains both positive and negative entries, the graded parts of the module have infinite dimension, and hence the Hilbert-Poincar\'e series is not well defined. That is why the proof of Theorem~\ref{mainTheo} cannot be generalized. Nevertheless, computational evidence suggests that the result is also true in this case. There are some infinite families of polynomials (such as the free ones, and those of the form $F=G\cdot H$ where $G$ satisfies the theorem and $H$ defines a hyperplane arrangement) for which alternative proofs exist. This justifies the following conjecture.
 
\begin{conjet} Theorem~\ref{mainTheo} also holds for aribtrary weight vector $\u\in \bz^n$. \end{conjet}

\section{Towards the General Case}\label{generalCase}

When $f\in S$ is not quasi-homogeneous, the module of logarithmic derivations
is not graded and thus the exponents and the graded betti numbers of $D(f)$ are not
well defined. In this section we show what can be expected in this general setting.
Let us start with an example in dimension 3.

\begin{ex}\label{lastEx}
Consider $f = x^2 z + y^3 + z^4 \in \K[x,y,z] = S$. A free resolution of $D(f)$ can be
efficiently computed with, for instance, the computer algebra system {\sc Singular} \cite{GPS09}.
\begin{equation}\label{exRes}
  0 \longleftarrow D(f) \stackrel{\varphi_0}{\longleftarrow}
  S(-1) \oplus S(-2) \oplus S(-3) \oplus S(-3)
  \stackrel{\varphi_1}{\longleftarrow} S(-5) \longleftarrow 0
\end{equation}

Although $D(f)$ is not graded and therefore the resolution is not homogeneous,
the notion of degree makes sense. Thus we write
$\varphi_0: S(-1) \oplus S(-2) \oplus S(-3) \oplus S(-3) \longrightarrow D(f)$ to emphasize
that $D(f)$ is generated by four elements of degrees $1,2,3,3$. In other words,
$\varphi_i$ is compatible with the usual filtration given by the degree of
a polynomial, that is, $\varphi_i(m)$ has degree $\leq k$ if $m$ has degree
$\leq k$, in the corresponding free graded $S$-module.

Note that the alternating sum of the degrees that appear in the above
resolution is the degree of $f$. However, replacing the basis
$\{ e_1,e_2,e_3,e_4 \}$ of the first free module by $\{e_1+e_2,e_2,e_3,e_4\}$,
one obtains a new free resolution of $D(f)$ where the first free module is
$S(-2) \oplus S(-2) \oplus S(-3) \oplus S(-3)$ and the rest being the same.
Now the alternating sum is different from the degree of the polynomial.
\end{ex}

The previous example tells us that Theorem \ref{mainTheo} can not be
stated for any resolution of $D(f)$ when $f$ is not quasi-homogeneous.
We shall show that it is always possible to find a resolution such that
the alternating sum of some ``exponents'' equals the degree of the polynomial
and also we will see how to compute such a resolution starting from a system
of generators of $D(f)$.

\subsection{Some Well-known Results on Homogenization}

Let us denote by $S^h = \K[x_1,\ldots,x_n,h] = S[h]$ the ring of polynomials
in $n+1$ variables. Here we consider on $S^h$ the classical graded structure given
by the degree of a polynomial, see Example \ref{ejemHom}. On $(S^h)^k$ the graded
structure is the one induced by $S^h$ considering that the elements in the canonical basis
have degree zero.

Given $M$ a $S$-submodule of $S^k$ and $m\in M$, we denote by $m^h \in (S^h)^k$
the homogenization of $m$ with respect to $h$. The homogenization of a
module is given by $M^h := S^h \langle m^h \mid m\in M \rangle$. Note that
if $\{m_1,\ldots,m_r\}$ is a system of generators of $M$, then
$\langle m_1^h,\ldots,m_r^h \rangle$ is in general strictly contained in $M^h$.

\begin{lemma}
The following properties about homogenization are satisfied.
\begin{enumerate}
\item $(m^h)_{|h=1} = m$.
\item Given $\xi\in (S^h)^k$ homogeneous, $\exists\l\in\mathbb{N}$ such that $h^\l (\xi_{|h=1})^h = \xi$.
\item $(M\cap N)^h = M^h \cap N^h$.
\item Let $\{m_1,\ldots,m_r\}$ be a Gr\"obner basis of $M$ with respect to
a graded monomial ordering. Then $M^h = S^h \langle m_1^h,\ldots,m_r^h\rangle$.
\end{enumerate}
\end{lemma}

In the sequel we will make use of this lemma without an explicit reference.

\subsection{Computing $\X(D(f)^h)$}

\begin{lemma}\label{lemmaDfh}
$D(f^h) \cap S^h \langle \partial_{x_1},\ldots,\partial_{x_n} \rangle = D(f)^h$.
\end{lemma}

\begin{proof}
Let $f=f_1^{e_1}\cdots f_r^{e_r}$ be the decomposition of $f$ into irreducible factors.
Let us denote $F=f^h$ and $F_i=f_i^h$. Consider
$\delta = \sum_{i=1}^n A_i \partial_{x_i} \in D(f^h)$ with $A_i\in S^h$ homogeneous of
the same degree. There exist $G_1,\ldots,G_r\in S^h$ such that
$\sum_{i=1}^n A_i \frac{\partial F_j}{\partial x_i} = G_j \cdot F_j^{e_j}$, for
all $j=1,\ldots,r$. Substituting $h=1$, one deduces that $\delta_{|h=1}$ belongs to
$D(f)$ and hence $\delta = h^\l (\delta_{|h=1})^h \in D(f)^h$.

Now take $\delta\in D(f)$. Then $\delta(f_j) = g_j\cdot f_j^{e_j}$ for
some $g_j\in S$, $j=1,\ldots,r$. We have that $\delta^h (F_i) =
h^{\l_i} \delta(f_i)^h = h^{\l_i} g_i^h F_i^{e_i} \in \langle F_i^{e_i}\rangle$.
This means that $\delta^h\in D(F)$ and the proof is complete.
\end{proof}

\begin{theo}\label{notMainTheo}
$\X(D(f)^h) = d$.
\end{theo}

\begin{proof}
Let us first assume that $f=f_1^{e_1}$ is $x_1$-general of order
$d$, the degree of the polynomial. Thus $F=f^h$ is also $x_1$-general.
From Theorem \ref{XDQhom}, $\X(D(F))=d$. The invariant $\X$ of
$S^h\langle \partial_{x_1},\ldots,\partial_{x_n}\rangle$ is zero,
since the module is isomorphic to $S^h(0)^n$. We have already seen in
the proof of Theorem \ref{XDQhom} that $\X\big(Der_\K(S)/(D(F)+S^h\langle \partial_{x_1},\ldots,\partial_{x_n}\rangle)\big) = 0$.
Hence we have,
\begin{equation*}
\begin{split}
  \X(D(f)^h)  =\, & \X(D(f^h)\cap S^h\langle \partial_{x_1},\ldots,\partial_{x_n}\rangle) =
  \X(D(F)) + \X(S^h \langle \partial_{x_1},\ldots,\partial_{x_n}\rangle) -\\
  & - \X(D(F) + S^h\langle \partial_{x_1},\ldots,\partial_{x_n}\rangle) = d + 0 - 0 = d.
\end{split}
\end{equation*}

Suppose the result is true for $f$ and $g$ with no common factors. Then
\begin{equation*}
\begin{split}
  \X(D(f g)^h) =\, & \X(D(f)^h\cap D(g)^h) = \X(D(f)^h) + \X(D(g)^h) -\\
  & - \X(D(f)^h + D(g)^h) = \deg(f) + \deg(g) - 0.
\end{split}
\end{equation*}
\end{proof}


\subsection{The Homogenization of a Resolution}

Let $M$ be an $S$-submodule of $D\simeq S^n$ and consider a free resolution
$0\leftarrow M \leftarrow F_\bullet \leftarrow 0$. As in Example \ref{lastEx},
each free module can be written as $F_p = \bigoplus_{i=1}^{r_p} S(-d_i^p)$ in
such a way that $\varphi_p: F_p \to F_{p-1}$ respects the filtration given by
the degree. Consider $F_p^h = \bigoplus_{i=1}^{r_p} S^h (-d_i^p)$ for $p\geq 0$,
$F_{-1}^h = M^h$ and define
$$
  \varphi_p^h : F_p^h \longrightarrow F_{p-1}^h:\
   e_i  \mapsto  \varphi_p^h(e_i):=\varphi_p(e_i)^h\,,
$$
where $\varphi_p(e_i)^h$ is the homogenization of $\varphi_p(e_i)$ in the
free module $F_{p-1}^h = \bigoplus_{i=1}^{r_{p-1}} S^h (-d_i^{p-1})$. Thus
one can obtain a new sequence of homogeneous morphisms of degree zero
$$
  0 \longleftarrow M^h \stackrel{\varphi_0^h}{\longleftarrow} F_0^h
  \stackrel{\varphi_1^h}{\longleftarrow} \cdots
  \stackrel{\varphi_\l^h}{\longleftarrow} F_\l^h \longleftarrow 0
$$
which a priori does not define a resolution of $M^h$.
\begin{prop}\label{homRes}
Using the above notation, we have that $(\ker \varphi_p)^h = \ker \varphi_p^h$
and $\im \varphi_p^h \subseteq (\im \varphi_p)^h$, that is, the homogenization
of the given resolution becomes a complex.
Moreover if $\im \varphi_p^h \supseteq (\im \varphi_p)^h$ holds, then this complex is
a homogeneous free resolution of $M^h$.
\end{prop}

\begin{proof}
Consider $a\in\ker \varphi_p$, then $\varphi_p(a)=0$. Thus $\varphi_p^h(a^h)$
is a homogeneous element in $F_{p-1}^h$ which is zero after making the substitution $h=1$ and therefore
is zero. Conversely, let $\xi \in F_{p}^h$ a homogeneous element such that $\varphi_p^h(\xi)=0$.
Then $\varphi_p(\xi_{|h=1}) = \varphi_p^h(\xi)_{|h=1} = 0$. Hence $\xi = h^\l (\xi_{|h=1})^h \in (\ker \varphi_p)^h$.

The second part of the statement is obvious.
\end{proof}

As a conclusion of the previous discussion we can algorithmically obtain a resolution of $D(f)$
which satisfies Theorem \ref{mainTheo} without any assumption on $f$. Compute
a resolution of $D(f)$ with respect to a degree monomial ordering and homogenized the resolution
in order to obtain a resolution of $D(f)^h$ following Proposition \ref{homRes}.
Finally, apply Theorem \ref{notMainTheo} and the additivity of $\X$.

\begin{ex}
Here we continue with Example \ref{lastEx}. The following is the matrix expression of the morphisms in
the resolution (\ref{exRes}).
\begin{small}
$$
\varphi_0 = \left(
\begin{array}{cccc}
9x & 3y^2 & 9z^3 & 0\\
8y & -2xz & -2xy & 4z^3+x^2\\
6z & 0 & -6xz & -3y^2
\end{array}\right);\quad
\varphi_1 = \left(
\begin{array}{c}
xy^2\\
-12z^3-3x^2\\
4y^2\\
-6xz
\end{array}\right)
$$
\end{small}

It can be checked that $D(f)^h = (\im \varphi_0)^h = \im \varphi_0^h$ and thus
the resolution verifies the conditions in Proposition \ref{homRes} and thus Theorem
\ref{mainTheo} holds.
However, replacing the basis
$\{ e_1,e_2,e_3,e_4 \}$ of the first free module by $\{e_1+e_2,e_2,e_3,e_4\}$,
one obtains a new free resolution of $D(f)$ whose matrix expression becomes
\begin{small}
$$
\psi_0 = \left(
\begin{array}{cccc}
9x+3y^2 & 3y^2 & 9z^3 & 0\\
8y-2xz & -2xz & -2xy & 4z^3+x^2\\
6z & 0 & -6xz & -3y^2
\end{array}\right);\quad
\psi_1 = \left(
\begin{array}{c}
-xy^2\\
xy^2+12z^3+3x^2\\
-4y^2\\
6xz
\end{array}\right).
$$
\end{small}

Now $\,\im \psi_0^h \nsupseteq(\im \psi_0)^h$ and Proposition~\ref{homRes} does not apply.  
\end{ex}

\section*{Acknowledgements}

We would like to thank Enrique Artal for
their constant support and motivation in our work over the years.
Also we wish to express our gratitude to Jos\'{e} Ignacio Cogolludo for his proofreading
and constant advise.


The authors are partially supported by the Spanish project MTM2007-67908-C02-01 and
``E15 Grupo Consolidado Geometr\'ia'' from the government of Arag\'on.
Also first author is partially supported by the ERC Starting Grant TGASS and
second author is partially supported by FQM-333.


\bibliography{./refbettiNumDQ}
\bibliographystyle{plain}

\end{document}